\magnification1080
\def\arXiv{oui}
  
 \def\oui{oui} 
  
\ifx\arXiv\oui
\else
 \pdfpagewidth=210truemm
 \pdfpageheight=297truemm 
\fi
  
  %
  %

  %
  %

  \catcode`@=12 

 \def\defrefnote#1{\definexref{#1}{{\the\footnotenumber}}{refnotes}}

  %
  %


\ifx\couleurs\oui
\input graphicx
 \pdfpagewidth=210truemm
 \pdfpageheight=297truemm 
 \voffset=-5mm
\fi

\input eplain.tex
\expandafter\def\expandafter\newdimen\expandafter{\newdimen}

\ifx\couleurs\oui
\beginpackages
\usepackage{color}
\endpackages 
 \pdfpagewidth=210truemm
 \pdfpageheight=297truemm 
\long\def\rge#1{{\color{red}#1}}

\definecolor{bleu-iecn}{cmyk}{.98,.13,.1,.55}

\else
\long\def\rge#1{#1}

\fi

\makeatletter
\def\numberedfootnote{%
ÊÊ\global\advance\footnotenumber by 1
ÊÊ\@eplainfootnote{{\number\footnotenumber}}%
}%
\def\makecolumns#1/#2 {\par \begingroup
ÊÊ \@columndepth = #1
ÊÊ \advance\@columndepth by -1
ÊÊ \divide \@columndepth by #2
ÊÊ \advance\@columndepth by 1
ÊÊ \@linestogoincolumn = \@columndepth
ÊÊ \@linestogo = #1
ÊÊ \currentcolumn = 1
ÊÊ \def\@endcolumnactions{%
ÊÊÊÊÊÊ\ifnum \@linestogo<2
ÊÊÊÊÊÊÊÊ \the\crtok \egroup \endgroup \par 
ÊÊÊÊÊÊ\else
ÊÊÊÊÊÊÊÊ \global\advance\@linestogo by -1
ÊÊÊÊÊÊÊÊ \ifnum\@linestogoincolumn<2
ÊÊÊÊÊÊÊÊÊÊÊÊ\global\advance\currentcolumn by 1
ÊÊÊÊÊÊÊÊÊÊÊÊ\global\@linestogoincolumn = \@columndepth
ÊÊÊÊÊÊÊÊÊÊÊÊ\the\crtok
ÊÊÊÊÊÊÊÊ \else
ÊÊÊÊÊÊÊÊÊÊÊÊ&\global\advance\@linestogoincolumn by -1
ÊÊÊÊÊÊÊÊ \fi
ÊÊÊÊÊÊ\fi
ÊÊ }%
ÊÊ \makeactive\^^M
ÊÊ \letreturn \@endcolumnactions
ÊÊ \@columnwidth = \hsize
ÊÊÊÊ \advance\@columnwidth by -\parindent
ÊÊÊÊ \divide\@columnwidth by #2
ÊÊ \penalty\abovecolumnspenalty
ÊÊ \noindent 
ÊÊ \valign\bgroup
ÊÊÊÊ &\hbox to \@columnwidth{\strut \hsize = \@columnwidth ##\hfil}\cr
}%
\makeatother

\lefteqnumbers
   \def\testd{oui}
   \def\choixlat{\ifx\numadroite\testd\righteqnumbers
            \else  \lefteqnumbers\fi}
    \choixlat

\catcode`@=\letter
\def\@eplainfootnote#1{\let\@sf\empty 
  \ifhmode\edef\@sf{\spacefactor\the\spacefactor}\/\fi
  \global\advance\hlfootlabelnumber by 1
  \hlstart@impl{foot}{\hlfootlabel}%
  \hldest@impl{footback}{\hlfootbacklabel}%
  \hbox{$^{(#1)}$}%
  \hlend@impl{foot}%
  \@sf\vfootnote{#1.}%
}%
\catcode`@=\other

  \interfootnoteskip=0pt
  
  \everyfootnote={\eightpoint\leftskip=5truemm\rightskip5truemm}
  
  \hsize150truemm\vsize 240truemm\hoffset=5truemm

  \pretolerance=500\tolerance=1000\brokenpenalty=5000
  \parindent3mm
  
  \countdef\temps=170
  \temps=\time
  \countdef\nminutes=171{\nminutes = \time}
  \countdef\nheure=172
  \def\heure{\begingroup                     
     \temps = \time \divide\temps by 60
     \nheure = \temps                        
     \nminutes = \time
     \multiply\temps by 60
     \advance\nminutes by -\temps            
     \ifnum\nminutes<10 \toks1 = {0}%
     \else\toks1 = {}%
     \fi
     \number\nheure h\the\toks1 \number\nminutes  
  \endgroup}%

  \newcount\chstart
  \chstart=\pageno
 \headline={\ifnum\pageno=\chstart {\hfill} \else {\hss \tenrm --\ \folio\ --\hss}\fi}
  \footline={\hfill}
  \normalbaselines
  \frenchspacing
    \def\dater{\vglue-10mm\rightline{(\the\day/\the\month/\the\year)}}
  \def\dateheure{\vglue-10mm\rightline{(\the\day/\the\month/\the\year,\ \heure)}}

  \newif\ifpagetitre \pagetitretrue
\newtoks\hautpagetitre \hautpagetitre={\hfill}
\newtoks\baspagetitre \baspagetitre={\hfill}
\newtoks\auteurcourant \auteurcourant={\hfill}
\newtoks\titrecourant \titrecourant={\hfill}
\newtoks\hautpagegauche
\newtoks\hautpagedroite
\newtoks\hautpagemilieu
\hautpagemilieu={\tenrm\hfil -- \folio\ -- \hfil}
\hautpagegauche={\ifx\midfolio\oui\the\hautpagemilieu\else\tenrm\folio\hfill\the\auteurcourant\hfill\fi}
\hautpagedroite={\ifx\midfolio\oui\the\hautpagemilieu\else\hfill\the\titrecourant\hfill\tenrm\folio\fi}
\newtoks\baspagegauche \baspagegauche={\hfil}
\newtoks\baspagedroite \baspagedroite={\hfil}
\headline={\ifpagetitre\the\hautpagetitre
\else\ifodd\pageno\the\hautpagedroite\else\the\hautpagegauche\fi\fi }
\footline={\ifpagetitre\the\baspagetitre
\else\ifodd\pageno\the\baspagedroite
\else\the\baspagegauche\fi\fi \global\pagetitrefalse}

\def\pageblanche{\vfill\eject\pagetitretrue
\null\vfill\eject
\pagetitretrue
}
\def\chgtpage{\ifodd\pageno \else
\pageblanche \fi \pagetitretrue\titreun=0\footnotenumber=0}

\def\chgtpageincrtitreun{\ifodd\pageno \else
\pageblanche \fi \pagetitretrue\footnotenumber=0}

\def\majnombres{\ifodd\pageno \else
\pageblanche \fi \pagetitretrue\hautpoly\titreun=0\footnotenumber=0}

\def\hautspages#1#2{\auteurcourant={\ninepcap#1}\titrecourant={\nineit#2}}

\ifnum\chstart=\pageno \pagetitretrue\fi
  

  \def\PAR{\par}
  
  \def\leftnote#1{\vadjust{\setbox1=\vtop{\hsize 20mm\parindent=0pt\eightpoint
  \baselineskip=9pt\rightskip=4mm plus 4mm\vskip-4.7mm#1}\hbox{\kern-2cm\smash{\box1}}}}

  
  \def\raggedcenter{\leftskip=20pt plus 10em  
       \rightskip=\leftskip 
        \parfillskip=0pt 
         \spaceskip=.3333em \xspaceskip=.5em 
          \pretolerance=9999 \tolerance=9999
           \hyphenpenalty=9999 \exhyphenpenalty=9999 }
           
  \def\titrecentre#1{{\parindent0mm\raggedcenter
       \spaceskip=.6em plus .2em minus .2em\xspaceskip=.6em plus .2em minus .2em
        \tit#1\par}}
        


  \def\oui{oui}
  
\def\fontetitreun{\ifx\paradouze\oui\douzepts\gpdouze\twelvebf\textfont1=\twelveib\else
\quatorzepts\gpquatorze\fourteenbf\fi}

\def\fontetitreunl{\douzepts\textfont1=\twelveib\scriptfont1=\tenib\fourteenti}
 
 \def\fontetitredeux{\textfont1=\eleveni\ifx\paradouze\oui\onzepts\scriptfont1=\ninei\elevenit\else
                        \douzepts\twelveit\fi}
 
   \def\fontetitredeuxb{\ifx\paradouze\oui\onzepts\eleventi\gponze\textfont1=\elevenib\scriptfont1=\nineib
                         \else\douzepts\twelveti\scriptfont1=\twelveib\scriptfont1=\tenib\gpdouze\fi}
                         
\def\fontetitredeuxl{\onzepts\textfont1=\elevenbf\scriptfont1=\ninebf\twelvebf}
  
\def\fontetitretrois{\textfont0=\elevenrm\scriptfont0=\eightrm\textfont1=\eleveni
                      \scriptfont1=\eighti\scriptscriptfont1=\sixi\elevenit}
                      
\def\fontetitrequatre{\textfont0=\elevenrm\scriptfont0=\eightrm\textfont1=\eleveni
                      \scriptfont1=\eighti\scriptscriptfont1=\sixi\elevenrm}
  
  \newcount\titreun\titreun=0
  \newcount\titredeux\titredeux=0
  \newcount\titretrois\titretrois=0
  \newcount\titrequatre\titrequatre=0
  \newcount\enonce\enonce=0
  
  \def\incr#1{\global\advance#1 by 1 {\the #1}}
  \def\avance#1{\global\advance#1 by 1}
  \def\init#1{\global#1=0}
  
  \long\def\Indentation#1#2{\setbox10=\hbox{\fontetitreun#1}
  	                    \ifdim\wd10 < 4mm
                         \setbox10=\hbox to 4mm{\box10\hfill}
                       \else\ifdim\wd10 < 6mm
                         \setbox10=\hbox to 6mm{\box10\hfill}
  	                    \else\ifdim\wd10 < 8mm
                         \setbox10=\hbox to 8mm{\box10\hfill}
                       \else\ifdim\wd10 < 12mm
                         \setbox10=\hbox to 12mm{\box10\hfill}
                       \fi\fi\fi\fi
                       \dimen10=\hsize
                       \advance \dimen10 by -\wd10
                       \noindent \box10 %
                       \ignorespaces
                       \hbox{\vtop{\hsize=\dimen10\raggedright\noindent\fontetitreun#2}}}

  \long\def\paraun#1{\removelastskip\par\medskip\goodbreak\vskip0pt plus.01\vsize\penalty-100
                \vskip0pt plus-.01\vsize
  	              \init{\titredeux}\ifnum\optionparag=1{\init\eqnumber\init\enonce}\else{}\fi
                  \goodbreak{\fontetitreun
  	                \Indentation{\incr{\titreun}.\ }{\fontetitreun #1\par}}\nobreak\medskip}

 %
 %
 \long\def\paraunc#1{\removelastskip\par\bigskip\goodbreak\vskip0pt plus.01\vsize\penalty-100
                \vskip0pt plus-.01\vsize
  	              \init{\titredeux}
                 \ifnum\optionparag=1{\init{\eqnumber}\init\enonce}\else{}\fi
                  \goodbreak
  	                {\parindent0mm\raggedcenter\fontetitreun\incr{\titreun}.\ 
                     \fontetitreun #1\par}\nobreak\medskip}
                     
\newtoks\titreunl
\titreunl={\ifnum\titreun=1{I}\fi%
\ifnum\titreun=2{II}\fi%
\ifnum\titreun=3{III}\fi%
\ifnum\titreun=4{IV}\fi%
\ifnum\titreun=5{V}\fi%
\ifnum\titreun=6{VI}\fi%
\ifnum\titreun=7{VII}\fi%
\ifnum\titreun=8{VIII}\fi%
\ifnum\titreun=9{IX}\fi%
\ifnum\titreun=10{X}\fi%
\ifnum\titreun=11{XI}\fi%
\ifnum\titreun=12{XII}\fi%
\ifnum\titreun=13{XIII}\fi%
}
\long\def\paraunl#1{\removelastskip\par\bigskip\bigskip\goodbreak\vskip0pt plus.01\vsize\penalty-100
                \vskip0pt plus-.01\vsize
  	              \init{\titredeux}\ifnum\optionparag=1{\init\eqnumber\init\enonce}\else{}\fi
                  \goodbreak{\fontetitreunl
  	                \Indentation{\global\advance\titreun by 1{\the\titreunl}.\ }{\fontetitreunl #1\par}}\nobreak\smallskip}

  
  \long\def\paradeux#1{\init{\titretrois}\vskip0pt plus.01\vsize\penalty-10
                \vskip0pt plus-.01\vsize\ifx \elie\oui\medskip\ifnum\titredeux>0\medskip\fi\fi
                 \Indentation{\fontetitredeux\the\titreun${\cdot}$\incr{\titredeux}.}
                              {\fontetitredeux\textfont1=\eleveni#1}\nobreak\par }
  
  \long\def\paradeuxb#1{\init{\titretrois}\vskip0pt plus.001\vsize\penalty-10
                \vskip0pt plus-.01\vsize{\ifx \elie\oui\medskip\ifnum\titredeux>0\medskip\fi\fi
                  \Indentation
  {\fontetitredeuxb\the\titreun${\cdot}$\incr{\titredeux}.}{ \fontetitredeuxb#1}}\nobreak
\smallskip}

\newtoks\titredeuxl
\titredeuxl={\ifnum\titredeux=1{A}\fi%
\ifnum\titredeux=2{B}\fi%
\ifnum\titredeux=3{C}\fi%
\ifnum\titredeux=4{D}\fi%
\ifnum\titredeux=5{E}\fi%
\ifnum\titredeux=6{F}\fi%
\ifnum\titredeux=7{G}\fi%
\ifnum\titredeux=8{H}\fi%
\ifnum\titredeux=9{I}\fi%
\ifnum\titredeux=10{J}\fi%
\ifnum\titredeux=11{K}\fi%
\ifnum\titredeux=12{L}\fi%
\ifnum\titredeux=13{M}\fi%
}
 \long\def\paradeuxl#1{\init{\titretrois}\vskip0pt plus.001\vsize\penalty-10
                \vskip0pt plus-.01
                \vsize \bigskip%
                  \Indentation
     {\fontetitredeuxl\global\advance\titredeux by 1
  \quad \the\titreunl${\cdot}$\the\titredeuxl.}{ \fontetitredeuxl#1}
  \removelastskip\nobreak\smallskip}
  

  \long\def\paratrois#1{\init{\titrequatre}\ifdim\lastskip<\smallskipamount
                \removelastskip\smallskip\fi
                 \vskip0pt plus.01\vsize\penalty-10
                  \vskip0pt
plus-.01\vsize{\ifx \elie\oui\ifnum\titretrois>0\medskip\fi\fi
\Indentation{\fontetitretrois\the\titreun${\cdot}$\the\titredeux${\cdot}$\incr{\titretrois}.\ }
  {\hskip0mm\baselineskip=14pt\fontetitretrois#1}\nobreak\smallskip}}
  
  
  \long\def\paratroisl#1{\init{\titrequatre}\ifdim\lastskip<\smallskipamount
                \removelastskip\fi
                 \vskip0pt plus.01\vsize\penalty-10
                  \vskip0pt
plus-.01\vsize\ifx \elie\oui\bigskip
\fi
\Indentation{\fontetitretrois\quad \quad \the\titreunl{${\cdot}$}\the\titredeuxl${\cdot}$\incr{\titretrois}.\ }
  {\hskip0mm\fontetitretrois#1}\nobreak\smallskip}


  \long\def\paraquatre#1{\ifdim\lastskip<\smallskipamount
                \removelastskip\smallskip\fi
                 \vskip0pt plus.01\vsize\penalty-10
                  \vskip0pt
                  plus-.01\vsize\par
 
\Indentation{\fontetitrequatre \the\titreun{${\cdot}$}\the\titredeux${\cdot}$\the\titretrois${\cdot}$\incr{\titrequatre}.\ }
{\hskip0mm\fontetitrequatre#1}\nobreak\smallskip}


\newtoks\titrequatrel
\titrequatrel={\ifnum\titrequatre=1{a}\fi%
\ifnum\titrequatre=2{b}\fi%
\ifnum\titrequatre=3{c}\fi%
\ifnum\titrequatre=4{d}\fi%
\ifnum\titrequatre=5{e}\fi%
\ifnum\titrequatre=6{f}\fi%
\ifnum\titrequatre=7{g}\fi%
\ifnum\titrequatre=8{h}\fi%
\ifnum\titrequatre=9{i}\fi%
\ifnum\titrequatre=10{j}\fi%
\ifnum\titrequatre=11{k}\fi%
\ifnum\titrequatre=12{l}\fi%
\ifnum\titrequatre=13{m}\fi%
}
\long\def\paraquatrel#1{\ifdim\lastskip<\smallskipamount
                \removelastskip\smallskip\fi
                 \vskip0pt plus.01\vsize\penalty-10
                  \vskip0pt
                  plus-.01\vsize{\bigskip
\Indentation{\global\advance\titrequatre by 1
\fontetitrequatre\quad \quad \quad \the\titreunl${\cdot}$\the\titredeuxl${\cdot}$\the\titretrois${\cdot}$\the\titrequatrel.\ }
{\hskip0mm\fontetitrequatre#1}\nobreak\smallskip}}

\ifx\optionkeys\oui
\def\drefun#1{\definexref{¤#1}{{\the\titreun}}{}} 
\def\drefdeux#1{\definexref{¤#1}{{\the\titreun}.{\the\titredeux}}{}}
\def\dreftrois#1{\definexref{¤#1}{{\the\titreun}.{\the\titredeux}.{\the\titretrois}}{}}
\else
\def\drefun#1{\definexref{prg#1}{{\the\titreun}}{}} 
\def\drefdeux#1{\definexref{prg#1}{{\the\titreun}.{\the\titredeux}}{}}
\def\dreftrois#1{\definexref{prg#1}{{\the\titreun}.{\the\titredeux}.{\the\titretrois}}{}}
\fi

%


  \long\def\propdeux#1#2#3#4{%
       \avance{\enonce}
       \leavevmode\edef\temp{#2}%
         \ifx\temp\empty 
          \else
           \definexref{#2}{#1~{\the\titreun.\the\enonce}}{enonces}
            \definexref{s#2}{{\the\titreun.\the\enonce}}{enonces}
             \fi
\smallskip
      \noindent{\bf#1\ {\bf\the\titreun.\the\enonce{#3}.}\enspace}{\sl#4\par}%
      \ifdim\lastskip<\medskipamount \removelastskip\penalty55\par \fi
   }

  \long\def\propun#1#2#3#4{%
      \avance{\enonce}
       \leavevmode\edef\temp{#2}%
        \ifx\temp\empty 
          \else
           \definexref{#2}{#1~{\the\enonce}}{enonces}
            \definexref{{s#2}}{{\the\enonce}}{enonces}
             \fi
   \par 
     \noindent{\bf#1\ {\bf\the\enonce{#3}.}\enspace}{\sl#4\par}%
     \ifdim\lastskip<\medskipamount \removelastskip\penalty55\medskip\fi
  }
  
  \long\def\prop#1#2#3#4{\ifnum\optionparag=1
                          \propdeux{#1}{#2}{\textfont1=\elevenib#3}{#4} \else\propun{#1}{#2}{\textfont1=\elevenib#3}{#4}\fi}

  \long\def\propt#1#2#3{\ifx\tpf\oui \prop{Th\'eo\-r\`eme}{#1}{#2}{#3}\par
                       \else\prop{Theorem}{#1}{#2}{#3}\par\fi}
  \long\def\Propt#1#2{\propt{#1}{}{#2}}
  \long\def\propl#1#2#3{\ifx\tpf\oui\prop{Lem\-me}{#1}{#2}{#3}\par
                         \else\prop{Lemma}{#1}{#2}{#3}\par\fi}
  \long\def\Propl#1#2{\propl{#1}{}{#2}}
  \long\def\propc#1#2#3{\ifx\tpf\oui\prop{Corol\-laire}{#1}{#2}{#3}\par
                         \else\prop{Corollary}{#1}{#2}{#3}\par\fi}
  \long\def\Propc#1#2{\propc{#1}{}{#2}}

  \long\def\propd#1#2#3{\ifx\tpf\oui\prop{D\'efi\-nition}{#1}{#2}{#3}\par
                       \else\prop{Definition}{#1}{#2}{#3}\par\fi} 
  
  \long\def\proptd#1#2#3{\ifx\tpf\oui\prop{Th\'eor\`eme et d\'efi\-nition}{#1}{#2}{#3}\par
                       \else\prop{Theorem and definition}{#1}{#2}{#3}\par\fi}


  
  \newcount\optionparag\optionparag=1
  
  \long\def\section#1#2{\ifnum\optionparag=1 \paraun{#2} 
                        \else\goodbreak{\fontetitreun
  	                \Indentation{#1.\ }{#2}}\nobreak\smallskip\fi}

  \def\eqconstruct#1{\ifnum\optionparag=1{\the\titreun\hbox{$\cdot$}#1}\else{#1}\fi}

  
  
  \def\numref{oui}  
  
  \newcount\mesref\mesref=0 
  \def\defbib#1{\ifx\numref\oui\global\advance\mesref by 1 \definexref{#1}{{\the
                 \mesref}}{}\else\definexref{#1}{#1}{}\fi}
  \def\bibtem#1{\defbib{#1}\item{\citer{#1}}}
  \def\citer#1{[\ref{#1}]}
  \def\citeplus#1#2{[\ref{#1}; #2]}

  
  \font\seventeenmsa=msam10 at 17pt    
  \font\fourteenmsa=msam10 at 14pt
  \font\twelvemsa=msam10 at 12pt
  \font\tenmsa=msam10                 
  \font\ninemsa=msam10 at 9pt 
  \font\eightmsa=msam10 at 8pt 
  \font\sevenmsa=msam7 
  \font\sixmsa=msam10 at 6pt
  \font\fivemsa=msam5
  \newfam\msafam\textfont\msafam=\tenmsa\scriptfont\msafam=\sevenmsa\scriptscriptfont\msafam=\fivemsa
  
  \font\seventeenbb=msbm10 at 17pt     
  \font\fourteenbb=msbm10 at 14pt
  \font\twelvebb=msbm10 at 12pt
  \font\tenbb=msbm10                   
  \font\ninebb=msbm10 at 9pt 
  \font\eightbb=msbm10 at 8pt 
  \font\sevenbb=msbm7 
  \font\sixbb=msbm10 at 6pt
  \font\fivebb=msbm5 
  \newfam\bbfam\textfont\bbfam=\tenbb\scriptfont\bbfam=\sevenbb\scriptscriptfont\bbfam=\fivebb
  \def\bb{\fam\bbfam\tenbb}%

  \font\seventeenscaln=eusm10 at 17pt   
  \font\twelvescaln=eusm10 at 12pt
  \font\tenscaln=eusm10                
  \font\ninescaln=eusm10 scaled 900
  \font\eightscaln=eusm10 scaled 800
  \font\sevenscaln=eusm10 scaled 700
  \font\sixscaln=eusm10 scaled 600
   
  \newfam\scalnfam\textfont\scalnfam=\tenscaln\scriptfont\scalnfam=\sevenscaln\scriptscriptfont\scalnfam=\sixscaln
  \def\scaln{\fam\scalnfam\tenscaln}%
  \def\scal{\scaln}
  
  \font\tenscalb=eusb10                

  \font\sevenscalb=eusb10 scaled 700

  \newfam\scalbfam\textfont\scalbfam=\tenscalb\scriptfont\scalbfam=\sevenscalb
  %
  
  %
  %
  \font\fourteenrm=cmr12 scaled 1200
  \font\elevenrm=cmr10 at 11pt
  \font\twelverm=cmr12
  \font\ninerm=cmr9
  \font\eightrm=cmr8      
  \font\sevenrm=cmr7
  \font\sixrm=cmr6

  \font\seventeenpcap=cmcsc10 at 17pt
  \font\tenpcap=cmcsc10                        
  \font\ninepcap=cmcsc9
  \font\eightpcap=cmcsc8
  \font\sevenpcap=cmcsc10 scaled 700
  
  \newfam\pcapfam\textfont\pcapfam=\tenpcap\scriptfont\pcapfam=\sevenpcap
  \def\pcap{\fam\pcapfam\tenpcap}
  
  \font\seventeenrm=cmbx12 scaled 1400

  \font\fourteenbf=cmbx10 scaled 1400
  
  \font\twelvebf=cmbx12
  \font\elevenbf=cmbx10 at 11pt
  \font\ninebf=cmbx9  
  \font\eightbf=cmbx8
  \font\sixbf=cmbx6
  
  \font\tengot=eufm10                           
   
  \font\eightgot=eufm10 at 8truept 
  \font\sevengot=eufm7 
  \font\sixgot=eufm10 at 6 truept 
   
  \newfam\gotfam
  \textfont\gotfam=\tengot\scriptfont\gotfam=\sevengot\scriptscriptfont\gotfam=\sixgot
  %

  
  \def\tit{%
  \textfont0=\seventeenrm\scriptfont0=\tenrm\def\rm{\fam0\seventeenrm}%
  \textfont1=\seventeenib\scriptfont1=\twelveib%
  \textfont2=\seventeensy\scriptfont2=\twelvesy\scriptscriptfont2=\ninesy
  \textfont3=\seventeenex
  \textfont\itfam=\seventeenti
  \def\it{\fam\itfam\seventeenti}%
  \textfont\bbfam=\seventeenbb \scriptfont\bbfam=\twelvebb
  \def\bb{\fam\bbfam\seventeenbb}%
  \textfont\msafam=\seventeenmsa\scriptfont\msafam=\twelvemsa
  \textfont\scalnfam=\seventeenscaln
  \def\pcap{\fam\pcapfam\seventeenpcap}
  \normalbaselineskip=25pt\normalbaselines\rm}

  \font\seventeenti=cmbxti10 scaled 1680
  
  \font\fourteenti=cmbxti10 at 14pt
  
  \font\twelveti=cmbxti10 scaled 1200
  \font\eleventi=cmbxti10 at 11pt

  %
  %
  \font\twelveit=cmti12	
  \font\elevenit=cmti10 scaled 1100
  \font\nineit=cmti9
  \font\eightit=cmti8
  \font\sevenit=cmti7

  %
  %
 
 \font\seventeenib=cmmib10 scaled 1680
  \font\fourteenib=cmmib10 scaled 1400
  \font\twelveib=cmmib10 scaled 1200
  \font\elevenib=cmmib10 scaled 1100
  \font\tenib=cmmib10
\font\eightib=cmmib10 scaled 800
  \font\nineib=cmmib10 scaled 900
\font\sevenib=cmmib10 scaled 700
\font\sixib=cmmib10 scaled 600
\font\fiveib=cmmib10 scaled 500

\ifx\ITAN\oui
\else
\innernewfam\cmmibfam
\textfont\cmmibfam=\tenib
\scriptfont\cmmibfam=\sevenib
\scriptscriptfont\cmmibfam=\fiveib
\def\ib{\fam\cmmibfam\tenib}
\fi

  %
  %
  \font\twelvei=cmmi10 scaled 1200
  \font\eleveni=cmmi10 scaled 1100
  \font\ninei=cmmi9
  \font\eighti=cmmi8 
  \font\seveni=cmmi7 			                
  \font\sixi=cmmi6
  
  \font\ninesl=cmsl9                    
  \font\eightsl=cmsl8 
  \font\sevensl=cmsl10 at 7pt

  \font\ninett=cmtt9                    
  \font\eighttt=cmtt8
  \font\seventt=cmtt10 scaled 700

  \font\seventeensy=cmsy10 scaled 1680    
  \font\fourteensy=cmsy10 scaled 1400
  \font\twelvesy=cmsy10 scaled 1176
  
  \font\ninesy=cmsy9                      
  \font\eightsy=cmsy8
  \font\sixsy=cmsy6
  \font\seventeenex=cmex10 at 17pt
  \font\fourteenex=cmex10 at 14pt
  \font\twelveex=cmex10 at 12pt
  \font\nineex=cmex10 at 9pt
  \font\eightex=cmex10 at 8pt
  \font\sevenex=cmex10 at 7pt
  \font\sixex=cmex10 at 6pt
  \font\fiveex=cmex10 at 5pt
  
   
  \font\fourteengp=cmmi10 at 14pt
  
  \font\twelvegp=cmmib10 at 12pt
  \font\elevengp=cmmib10 at 11pt
  \font\tengp=cmmib10                          
  \font\ninegp=cmmib10 at 9pt 
  \font\eightgp=cmmib8 
   
  \font\sixgp=cmmib6


  \def\gponze{\textfont0=\elevenbf\scriptfont0=\eightbf\scriptscriptfont0=\sixbf
           \textfont1=\elevengp\scriptfont1=\eightgp\scriptscriptfont1=\sixgp}
  \def\gpdouze{\textfont0=\twelvebf\scriptfont0=\tenbf\scriptscriptfont0=\ninebf
           \textfont1=\twelvegp\scriptfont1=\tengp\scriptscriptfont1=\ninegp}        
  
 \def\gpquatorze{\textfont0=\fourteenbf\scriptfont0=\twelvebf\scriptscriptfont0=\elevenbf
           \textfont1=\fourteengp\scriptfont1=\twelvegp\scriptscriptfont1=\elevengp}

  
  \expandafter\chardef\csname pre amssym.def at\endcsname=\the\catcode`\@
  \catcode`\@=11
  \def\undefine#1{\let#1\undefined}
  \def\newsymbol#1#2#3#4#5{\let\next@\relax
   \ifnum#2=\@ne\let\next@\msafam@\else
   \ifnum#2=\tw@\let\next@\bbfam@\fi\fi
   \mathchardef#1="#3\next@#4#5}
  \def\mathhexbox@#1#2#3{\relax
   \ifmmode\mathpalette{}{\m@th\mathchar"#1#2#3}%
   \else\leavevmode\hbox{$\m@th\mathchar"#1#2#3$}\fi}
  \def\hexnumber@#1{\ifcase#1 0\or 1\or 2\or 3\or 4\or 5\or 6\or 7\or 8\or
   9\or A\or B\or C\or D\or E\or F\fi}
  
  \def\setboxz@h{\setbox\z@\hbox}
  \def\wdz@{\wd\z@}
  \def\boxz@{\box\z@}
  
  \edef\msafam@{\hexnumber@\msafam}
  \mathchardef\dabar@"0\msafam@39
  
  \edef\bbfam@{\hexnumber@\bbfam}
  \def\widehat#1{\setboxz@h{$\m@th#1$}%
   \ifdim\wdz@>\tw@ em\mathaccent"0\bbfam@5B{#1}%
   \else\mathaccent"0362{#1}\fi}
  \def\widetilde#1{\setboxz@h{$\m@th#1$}%
   \ifdim\wdz@>\tw@ em\mathaccent"0\bbfam@5D{#1}%
   \else\mathaccent"0365{#1}\fi}
  \newsymbol\leqq 1335          
  \newsymbol\leqslant 1336
  \newsymbol\lessgtr 1337       
  \newsymbol\backprime 1038     
  \newsymbol\risingdotseq 133A  
  \newsymbol\fallingdotseq 133B 
  \newsymbol\succcurlyeq 133C   
  \newsymbol\geqq 133D          
  \newsymbol\geqslant 133E
  \newsymbol\nmid 232D
  \newsymbol\nexists 2040
  \newsymbol\smallsetminus 2272
  \newsymbol\varnothing 203F 
  \catcode`\@=\active

  \catcode`\@=11
  \newcount\typofr\typofr=1
  
  \catcode`\;=\active
  \def;{\ifnum\typofr=1\relax\ifhmode\ifdim\lastskip>\z@\unskip\fi
     \kern.2em\fi\string;\else\string;\fi}
  
  \catcode`\:=\active
  \def:{\ifnum\typofr=1\relax\ifhmode\ifdim\lastskip>\z@\unskip\fi
  \penalty\@M\ \fi\string:\else\string:\fi}
  
  \catcode`\!=\active
  \def!{\ifnum\typofr=1\relax\ifhmode\ifdim\lastskip>\z@\unskip\fi
     \kern.2em\fi\string!\else\string!\fi}
  
  \catcode`\?=\active
  \def?{\ifnum\typofr=1\relax\ifhmode\ifdim\lastskip>\z@\unskip\fi
     \kern.2em\fi\string?\else\string?\fi}

  \def\francais{\typofr=1\def\tpf{oui}}
  \def\anglais{\typofr=2\def\tpf{non}\def\english{oui}}
  \def\oui{oui}
  \francais
  
  \catcode`\@=12
  

\ifx\textures\oui
\def\raye #1|{\leavevmode\setbox1=\hbox{#1}%
\raise .5pt\hbox to \wd1{\xleaders\hbox{\rge{$ \sct / $}%
\kern 1pt}\hfill\kern -1pt }\kern -\wd1 \unhbox1\relax }

\def\barre #1|{\leavevmode\setbox1=\hbox{#1}%
\rlap{\Red\vrule height 2.4pt depth -1.2pt width \wd1}\Black \unhbox1\relax}
\else
\def\raye #1|{\leavevmode\setbox1=\hbox{#1}%
\raise .5pt\hbox to \wd1{\xleaders\hbox{\rge{$ \sct / $}%
\kern 1pt}\hfill\kern -1pt }\kern -\wd1 \unhbox1\relax }

\def\barre #1|{\leavevmode\setbox1=\hbox{#1}%
\rlap{\color{red}\vrule height 2.4pt depth -1.2pt width \wd1}\color{black} \unhbox1\relax}

\fi
  

  
  \def\og{\leavevmode\raise.24ex\hbox{$\scriptscriptstyle\langle\!\langle\>$}}    
  \def\fg{\leavevmode\raise.24ex\hbox{$\scriptscriptstyle\>\rangle\!\rangle$}}    

  \def\d{\,{\rm d}}

  \def\r{{\bb R}}
  
  \def\N{{\bb N}}

  \def\A{{\scal A}}

  \def\HH{{\scal H}}

  \def\M{{\scal M}}
  
  \def\O{{\scal O}}
  \def\P{{\scaln P}}

  \def\frac#1#2{{#1\over #2}}
  \def\di#1#2{\sct#1\atop{\sct#2}}

  \def\numero{n$^{\rm o}\thinspace$}

  \def\qedbox{$\rlap{$\sqcap$}\sqcup$}           
  \def\qed{\nobreak\hfill\penalty250 \hbox{}\nobreak\hfill\qedbox\par }

  \def\sumast{\mathop{{\sum}^*}}
  
  \def\numero{n$^{\rm o}\thinspace$}

  \def\pnu{p^\nu}

  \def\¤{\S\thinspace}

  \def\¥{$\bullet$ }
  
  
  \def\e{{\rm e}}

  \def\epsilon{\varepsilon}

  \def\phi{\varphi}
  \def\theta{\vartheta}
  \def\rho{\varrho}
  \def\dm{{\textstyle{1\over 2}}}
  \def\txt{\textstyle}
  \def\dsp{\displaystyle}
  \def\sct{\scriptstyle}
  \def\pf{\noi{\it Proof. }}
  \def\nid{\ifnum\typofr=1\par\noindent{\it D\'emonstration. }\else\pf\fi}
  \def\noi{\noindent}
  \def\rem{\ifnum\typofr=1\noi{\it Remarque.}\ \else\noi{\it Remark.}\ \fi}
  \def\rems{\ifnum\typofr=1\noi{\it Remarques.}\ \else\noi{\it Remarks.}\ \fi}

  \def\sset{\smallsetminus}

  \def\pp{{\rm pp}}
  
  \def\meas{\mathop{\rm meas}\nolimits}
  
  \def\1{{\bf 1}}
  \def\|{\Vert}

  \def\leq{\leqslant}
  \def\geq{\geqslant}

  \def\eg{{e.g.}}
  \newsymbol\subsetneqq 2324
  \newsymbol\subsetneq 2328

  \def\fl#1{\left\lfloor #1 \right\rfloor}

  \def\li{\mathop{\rm li}\nolimits}

  \def\log{\mathop{\rm log}\nolimits}
  \def\ft#1#2{{\txt{#1\over #2}}}




  \def\pmb#1{\setbox0=\hbox{#1}%
  \kern-.025em\copy0\kern-\wd0\kern.05em\copy0\kern-\wd0\kern-.025em\raise .0433em\box0 }

  
  \skewchar\eighti='177 \skewchar\sixi='177
  \skewchar\eightsy='60 \skewchar\sixsy='60
  
  \def\eightpoint{%
  \textfont0=\eightrm\scriptfont0=\sixrm\scriptscriptfont0=\fiverm
  \def\rm{\fam0\eightrm}%
  \textfont1=\eighti\scriptfont1=\sixi
  \scriptscriptfont1=\fivei\def\oldstyle{\fam1\seveni}%
  \textfont2=\eightsy\scriptfont2=\sixsy\scriptscriptfont2=\fivesy
  \textfont3=\eightex\scriptfont3=\sixex
  \textfont\itfam=\eightit
  \def\it{\fam\itfam\eightit}%
  \textfont\slfam=\eightsl
  \def\sl{\fam\slfam\eightsl}%
  \textfont\bbfam=\eightbb \scriptfont\bbfam=\sixbb\scriptscriptfont\bbfam=\fivebb
  \def\bb{\fam\bbfam\eightbb}%
  \textfont\msafam=\eightmsa\scriptfont\msafam=\sixmsa
  \textfont\scalnfam=\eightscaln
  \def\scaln{\fam\scalnfam\eightscaln}
  \textfont\ttfam=\eighttt
  \def\tt{\fam\ttfam\eighttt}%
\textfont\gotfam=\eightgot
  \textfont\bffam=\eightbf\scriptfont\bffam=\sixbf\scriptscriptfont\bffam=\fivebf
  \def\bf{\fam\bffam\eightbf}%
  \ifx\ITAN\oui\else\textfont\cmmibfam=\eightib
       \scriptfont\cmmibfam=\sixib
        \scriptscriptfont\cmmibfam=\fiveib
         \def\ib{\fam\cmmibfam\eightib}
   \fi
  \textfont\pcapfam=\eightpcap
  \def\pcap{\fam\pcapfam\eightpcap}
  \abovedisplayskip=2pt plus2pt minus 2pt
  \belowdisplayskip=2pt plus1pt minus 2pt
  \abovedisplayshortskip= 1pt plus 2pt minus 1pt
  \belowdisplayshortskip= 1pt plus 2pt minus 1pt
  \smallskipamount=2pt plus 1pt minus 2pt
  \medskipamount=3pt plus 2pt minus 2pt
  \bigskipamount=7pt plus 3pt minus 3pt
  \setbox\strutbox=\hbox{\vrule height 5pt depth 2pt width 0pt}%
  \normalbaselineskip=9pt\normalbaselines\rm}

  \def\({\left(}
  \def\){\right)}
  
  \def\footnoterule{\kern -2pt\hrule width 7truecm\kern 2.4pt}
  
  \def\xnotedef#1{\definexref{#1}{\noexpand\number\footnotenumber}{Note}}%

  
  
  \def\ninepoint{%
  \textfont0=\ninerm\scriptfont0=\sixrm\scriptscriptfont0=\fiverm
  \def\rm{\fam0\ninerm}%
  \textfont1=\ninei\scriptfont1=\sixi
  \scriptscriptfont1=\fivei\def\oldstyle{\fam1\ninei}%
  \textfont2=\ninesy\scriptfont2=\sixsy\scriptscriptfont2=\fivesy
  \textfont3=\nineex\scriptfont3=\sixex
  \textfont\itfam=\nineit
  \def\it{\fam\itfam\nineit}%
  \textfont\slfam=\ninesl
  \def\sl{\fam\slfam\ninesl}%
  \textfont\bbfam=\ninebb\scriptfont\bbfam=\sixbb\scriptscriptfont\bbfam=\fivebb
  \def\bb{\fam\bbfam\ninebb}%
  \textfont\msafam=\ninemsa\scriptfont\msafam=\sixmsa\scriptscriptfont\msafam=\fivemsa
  \textfont\scalnfam=\ninescaln
  \def\scaln{\fam\scalnfam\ninescaln}
  \textfont\ttfam=\ninett
  \def\tt{\fam\ttfam\ninett}%
  \textfont\bffam=\ninebf\scriptfont\bffam=\sixbf\scriptscriptfont\bffam=\fivebf
  \def\bf{\fam\bffam\ninebf}%
  \abovedisplayskip=3pt plus2pt minus 2pt
  \belowdisplayskip=3pt plus1pt minus 2pt
  \abovedisplayshortskip= 2pt plus 2pt minus 1pt
  \belowdisplayshortskip= 2pt plus 2pt minus 1pt
  \smallskipamount=2pt plus 1pt minus 2pt
  \medskipamount=3pt plus 2pt minus 2pt
  \bigskipamount=7pt plus 3pt minus 3pt
  \setbox\strutbox=\hbox{\vrule height 5pt depth 2pt width 0pt}%
  \normalbaselineskip=11pt plus.3pt minus.2pt\normalbaselines\rm}

  \def\sevenpoint{%
  \textfont0=\sevenrm\scriptfont0=\sixrm\scriptscriptfont0=\fiverm
  \def\rm{\fam0\sevenrm}%
  \textfont1=\seveni\scriptfont1=\sixi
  \scriptscriptfont1=\fivei\def\oldstyle{\fam1\seveni}%
  \textfont2=\sevensy\scriptfont2=\sixsy\scriptscriptfont2=\fivesy
  \textfont3=\sevenex\scriptfont3=\fiveex
  \textfont\itfam=\sevenit
  \def\it{\fam\itfam\sevenit}%
  \textfont\slfam=\sevensl
  \def\sl{\fam\slfam\sevensl}%
  \textfont\bbfam=\sevenbb \scriptfont\bbfam=\sixbb\scriptscriptfont\bbfam=\fivebb
  \def\bb{\fam\bbfam\sevenbb}%
  \textfont\msafam=\sevenmsa\scriptfont\msafam=\sixmsa
  \textfont\scalnfam=\sevenscaln
  \def\scaln{\fam\scalnfam\sevenscaln}
  \textfont\bffam=\sevenbf\scriptfont\bffam=\sixbf\scriptscriptfont\bffam=\fivebf
  \def\bf{\fam\bffam\sevenbf}%
  \textfont\ttfam=\seventt
  \abovedisplayskip=2pt plus2pt minus 2pt
  \belowdisplayskip=2pt plus1pt minus 2pt
  \abovedisplayshortskip= 1pt plus 2pt minus 1pt
  \belowdisplayshortskip= 1pt plus 2pt minus 1pt
  \smallskipamount=2pt plus 1pt minus 2pt
  \medskipamount=3pt plus 2pt minus 2pt
  \bigskipamount=7pt plus 3pt minus 3pt
  \setbox\strutbox=\hbox{\vrule height 5pt depth 2pt width 0pt}%
  \normalbaselineskip=9pt\normalbaselines\rm}

 \def\onzepts{%
 \textfont0=\elevenrm\scriptfont0=\ninerm
 \textfont1=\eleveni\scriptfont1=\ninei
}

\def\douzepts{%
  \textfont0=\twelverm\scriptfont0=\tenrm\def\rm{\fam0\twelverm}%
  \textfont1=\twelvei\scriptfont1=\teni%
  \textfont2=\twelvesy\scriptfont2=\tensy\scriptscriptfont2=\eightsy
  \textfont3=\twelveex
  \textfont\itfam=\twelveti
  \def\it{\fam\itfam\twelveti}%
  \textfont\bffam=\twelvebf\scriptfont\bffam=\tenbf\scriptscriptfont\bffam=\eightbf
  \def\bf{\fam\bffam\twelvebf}%
  \textfont\bbfam=\twelvebb \scriptfont\bbfam=\tenbb
  \def\bb{\fam\bbfam\twelvebb}%
  \textfont\msafam=\twelvemsa\scriptfont\msafam=\tenmsa
  \textfont\scalnfam=\twelvescaln
  \normalbaselineskip=15pt\normalbaselines\rm}

\def\quatorzepts{%
  \textfont0=\fourteenrm\scriptfont0=\twelverm\def\rm{\fam0\fourteenrm}%
  \textfont1=\fourteenib\scriptfont1=\twelveib%
  \textfont2=\fourteensy\scriptfont2=\twelvesy\scriptscriptfont2=\tensy
  \textfont3=\fourteenex
  \textfont\itfam=\fourteenti
  \def\it{\fam\itfam\fourteenti}%
  \textfont\bffam=\fourteenbf\scriptfont\bffam=\twelvebf\scriptscriptfont\bffam=\tenbf
  \def\bf{\fam\bffam\fourteenbf}%
  \textfont\bbfam=\fourteenbb \scriptfont\bbfam=\twelvebb
  \def\bb{\fam\bbfam\fourteenbb}%
  \textfont\msafam=\fourteenmsa\scriptfont\msafam=\twelvemsa
  \textfont\scalnfam=\twelvescaln
  \normalbaselineskip=18pt\normalbaselines\rm}


\def\AA{{\it Acta Arith.}}

\def\picture #1 by #2 (#3){\leavevmode\vbox to #2{
     \hrule width #1 height 0pt depth 0pt
      \vfill
       \special{picture #3}}}

\def\illustration #1 by #2 (#3) scaled #4{\dimen1=#2
  \divide\dimen1 by 1000
  \multiply\dimen1 by #4
  \vtop to \dimen1{\dimen1=#1
  \divide\dimen1 by 1000
  \multiply\dimen1 by #4
  \hsize=\dimen1\vss
  \noindent\special{illustration #3 scaled #4}}}

\ifx\couleurs\oui

\fi

\optionparag=2

\ifx\optionkeymacros\undefined\else \fi

\catcode`\Œ=\active\defŒ{{\aa}}       
\catcode`\º=\active\defº{\int}        
\catcode`\=\active\def{\c c}        
\catcode`\¶=\active\def¶{\partial}    
\catcode`\Ä=\active\defÄ{\oint}       
\catcode`\Æ=\active\defÆ{\triangle}   
\catcode`\Â=\active\defÂ{\neg}        
\catcode`\µ=\active\defµ{\mu}         
\catcode`\¿=\active\def¿{{\o}}        
\catcode`\¹=\active\def¹{\pi}         
\catcode`\Ï=\active\defÏ{{\oe}}       
\catcode`\§=\active\def§{{\ss}}       
\catcode`\ =\active\def {\dagger}     
\catcode`\Ã=\active\defÃ{\sqrt}       
\catcode`\·=\active\def·{\Sigma}      
\catcode`\Å=\active\defÅ{\approx}     
\catcode`\½=\active\def½{\Omega}      
\catcode`\£=\active\def£{{\it\$}}     
\catcode`\°=\active\def°{\infty}      
\catcode`\¤=\active\def¤{{\S}}        
\catcode`\¦=\active\def¦{{\P}}        
\catcode`\¥=\active\def¥{\bullet}     
\catcode`\»=\active\def»{\leavevmode\raise.585ex\hbox{\b a}}      
\catcode`\¼=\active\def¼{\leavevmode\raise.6ex\hbox{\b o}}        
\catcode`\­=\active\def­{\not=}       
\catcode`\²=\active\def²{\leq}        
\catcode`\³=\active\def³{\geq}        
\catcode`\Ö=\active\defÖ{\div}        
\catcode`\É=\active\defÉ{{\dots}}     
\catcode`\¾=\active\def¾{{\ae}}       
\catcode`\Ç=\active\defÇ{\og}         
\catcode`\Ò=\active\defÒ{``}          
\catcode`\Á=\active\defÁ{!`}          
\catcode`\¢=\active\def¢{\rlap/c}     
\catcode`\Ô=\active\defÔ{`}           
\catcode`\Õ=\active\defÕ{'}           


\catcode`\=\active\def{{\AA}}       
\catcode`\'=\active\def'{\c C}        
\catcode`\¯=\active\def¯{{\O}}        
\catcode`\¸=\active\def¸{\Pi}         
\catcode`\Î=\active\defÎ{{\OE}}       
\catcode`\®=\active\def®{{\AE}}       
\catcode`\×=\active\def×{\diamond}    
\catcode`\¡=\active\def¡{\accent'27}  
\catcode`\Ó=\active\defÓ{''}          
\catcode`\±=\active\def±{\pm}         
\catcode`\È=\active\defÈ{\fg}         
\catcode`\À=\active\defÀ{?`}          
\catcode`\Ð=\active\defÐ{--}          
\catcode`\Ñ=\active\defÑ{---}         


\catcode`\Š=\active\defŠ{\"a}        
\catcode`\'=\active\def'{\"e}        
\catcode`\•=\active\def•{\"{\i}}     
\catcode`\š=\active\defš{\"o}        
\catcode`\Ÿ=\active\defŸ{\"u}        
\catcode`\Ø=\active\defØ{\"y}        
\catcode`\å=\active\defå{\^A}        
\catcode`\€=\active\def€{\"A}        
\catcode`\…=\active\def…{\"O}        
\catcode`\†=\active\def†{\"U}        
\catcode`\‡=\active\def‡{\'a}        
\catcode`\Ž=\active\defŽ{\'e}        
\catcode`\'=\active\def'{\'{\i}}     
\catcode`\—=\active\def—{\'o}        
\catcode`\œ=\active\defœ{\'u}        
\catcode`\ƒ=\active\defƒ{\'E}        
\catcode`\æ=\active\defæ{\^E}        
\catcode`\é=\active\defé{\`E}        %
\catcode`\ˆ=\active\defˆ{\`a}        
\catcode`\=\active\def{\`e}        
\catcode`\"=\active\def"{\`{\i}}     
\catcode`\˜=\active\def˜{\`o}        
\catcode`\=\active\def{\`u}        
\catcode`\Ë=\active\defË{\`A}        
\catcode`\‹=\active\def‹{\~a}        
\catcode`\–=\active\def–{\~n}        
\catcode`\›=\active\def›{\~o}        
\catcode`\Ì=\active\defÌ{\~A}        
\catcode`\"=\active\def"{\~N}        
\catcode`\Í=\active\defÍ{\~O}        
\catcode`\‰=\active\def‰{\^a}        
\catcode`\=\active\def{\^e}        
\catcode`\"=\active\def"{\^{\i}}     
\catcode`\™=\active\def™{\^o}        
\catcode`\ž=\active\defž{\^u}        

\let\optionkeymacros\null

\vsize=255truemm
\voffset-5mm

\anglais

\def\pnu{p^\nu}
\optionparag=1
\def\paradouze{oui}

\def\auteur{RŽgis de la Bretche \& GŽrald Tenenbaum}
\def\titrart{Two upper bounds for the Erd\H os--Hooley Delta-function}
\hautspages{\auteur}{\titrart}
\dateheure
\titrecentre{\titrart}
\bigskip\medskip
\centerline {\auteur} 
\bigskip\bigskip
{\eightpoint\leftskip1cm\rightskip1cm
\noi{\bf Abstract.}  
For integer $n\geqslant 1$ and real $u$, let $\Delta(n,u):=|\{d:d\mid n,\,\e^u<d\leqslant \e^{u+1}\}|$. The Erd\H os--Hooley Delta-function is then defined by $\Delta(n):=\max_{u\in\r}\Delta(n,u).$ We improve the current upper bounds for the average and normal orders of this arithmetic function.\PAR
\medskip
\noi
{\bf Keywords:}  Concentration of divisors, average order, normal order, Waring's problem.\par
\smallskip 
\noi \bf 2020 Mathematics Subject Classification: \rm primary   11N37; secondary 11K65.\par }
\bigskip\medskip
\paraun{Introduction and statement of results}
For integer $n\geqslant 1$ and real $u$, put $$\Delta(n,u):=|\{d:d\mid n,\,\e^u<d\leqslant \e^{u+1}\}|,\quad\Delta(n):=\max_{u\in\r}\Delta(n,u).$$
 Introduced by Erd\H os \citer{Er73} (see also \citer{EN75}) and studied by Hooley \citer{Ho79}, the $\Delta$-function  proved very useful in several branches of number theory --- see, \eg, \citer{HT88} and \citer{Te86} for further references. If $\tau(n)$ denotes the total number of divisors of $n$, then $\Delta(n)/\tau(n)$ coincides with the concentration of the numbers $\log d$, $d|n$. In this work, we aim at improving the current upper bounds for the average and normal orders. In the former case, we consider weighted versions.\par 
For $A>0$, $y\geqslant 1$, $c>0$,  $\eta\in]0,1[$, we define
the class $\M_A(y,c,\eta)$ comprising those arithmetic functions
$g$ that are multiplicative, non-negative, and satisfy
the conditions
$$\leqalignno{&g(p^\nu)\leqslant A^\nu\qquad
(\nu\geqslant 1)&\eqdef{majgpnu}\cr
&(\forall\varepsilon>0)\quad g(n)\ll_\varepsilon n^\varepsilon\quad(n\geqslant 1)&\eqdef{majgn}\cr
&\dsp\sum_{p\leqslant
x}{g(p)}=y\li(x)+O\big(x\e^{-c(\log
x)^\eta}\big)\quad(x\geqslant 2).&\eqdef{condg}\cr}
$$
Here and in the sequel we reserve the letter  $p$ to designate a prime number.
Note for the sake of further reference that a theorem of Shiu \citer{Sh80} implies 
$$\sum_{n\leqslant x}g(n)\ll x(\log x)^{y-1}\eqdef{Shiu}$$
for any $g$ in $\M_A(y,c,\eta).$
\par 
Regarding average values of the $\Delta$-function, we consider the weighted sum
$$S(x;g):=\sum_{n\leqslant x}g(n) \Delta(n).$$
Here and throughout we let $\log_k$ denote the $k$-fold iterated logarithm.
\vskip-3mm 
 \Propt{th1}{Let $A>0$, $y\geqslant 1$, $c>0$,  $\eta\in]0,1[$, $g\in  \M_A(y,c,\eta)$, $a>\sqrt{2} \log 2\approx 0.980258$. We have
 $$S(x;g)
\ll x (\log x)^{2y-2}\e^{Êa\sqrt{\log_2x}}\qquad (x\geqslant 3).\eqdef{majDg}$$}

We note that by a different approach Koukoulopoulos (private communication)  obtained  a similar estimate for  $g=\1$ with $a=2.1$

When $y\geqslant 1$, \ref{th1} provides a small improvement over known estimates, for instance \hbox{\citeplus{HT88}{th.~70}} stating that, for anyy $\varepsilon>0$,
$$S(x;\1)\ll x\e^{(1+\varepsilon)\sqrt{2\log_2 x\log_3x}}\qquad (x\to\infty).$$\par 
When $y\leqslant \dm$, the proof of \citeplus{HT88}{th. 64} may be readily 
generalized to any $g\in\M_A(y,c,\eta)$   to yield
$$S(x;g)
\ll x (\log x)^{ y-1}(\log_2x)^{\delta(y,1/2)},\eqdef{majDg-ypetit}$$
where $\delta(u,v):=1$ if $u=v$, and $:=0$ otherwise.\par 
For $y\geqslant 1+\dm\sqrt{2}$, we may adapt {\it mutatis mutandis} \citeplus{HT88}{th. 67}  and derive 
$$S(x;g)
\ll x (\log x)^{ 2y-2}(\log_2x)^{\delta(y,1+\sqrt{2}/2)}.\eqdef{majDg-ygrand}$$\par 
Finally, we note that the proof of \citeplus{HT88}{th. 71} may be extended to get, for any fixed  $y<1$, 
$$S(x;g)\ll x (\log x)^{y-1}\exp\Big\{ {4\log 3\over 1-y}(\log_3x)^{2}\Big\}.\eqdef{majDg-y<1}$$
We omit further details since the relevant approaches are straightforward.
\par \goodbreak

As put forward by Hooley \citer{Ho79}, among other applications, average bounds such as  \eqref{majDg} may be employed to count solutions of certain Diophantine equations. For given $k\in\N^*$ and positive integers $c_j$, $\ell_j$ $(0\leqslant j\leqslant k)$ with $\ell_0=\min \ell_j= 2$,  we consider as in \citer{Ro11} the number $N(x)$ of solutions $({\ib m}, {\ib n})=(m_0, \ldots, m_k,n_0, \ldots, n_k)\in \N^{2k+2}$ of the system
$$\sum_{0\leqslant j\leqslant k} c_j m_j^{\ell_j}=\sum_{0\leqslant j\leqslant k}c_j n_j^{\ell_j}\leqslant x,\ m_0\neq n_0,$$
and let $V(x)$ denote the number of integers $n\leqslant x$ that are representable in the form
$$n=\sum_{0\leqslant j\leqslant k} c_j n_j^{\ell_j}.$$
 The case $k=2$, $c_0=c_1=c_2=1$, $\ell_1=\ell_2=4$ has been studied by the second author~\citer{Te86}.

Applying \ref{th1} with $y=1$ and following the approach displayed in \citer{Ro11}, we derive the next   corollary.
\Propc{cor}{Assume $\sum_{1\leqslant j\leqslant k} 1/\ell_j=\ft 12$ and let $a>\sqrt{2}\log 2$. Then, as $x\to\infty$, we have  $$\leqalignno{N(x)&\ll x\e^{a\sqrt{\log_2 x}},&\eqdef{majN}\cr
V(x)&\gg x\e^{-a\sqrt{\log_2 x}}.&\eqdef{minV}\cr}$$
}
We omit the details of the proof,  since they are almost identical to those in \citer{Ro11}.

\medskip
We next turn our attention to the normal order of the $\Delta$-function. We employ the mention~$\pp$ to indicate that a formula holds on a sequence of natural density 1. Improving on estimates of  Maier \& Tenenbaum \citer{MT84}, \citer{MT09}, Ford, Green \& Koukoulopoulos \citer{FGK19} recently claimed
$$\Delta(n)>(\log_2n)^{\gamma_1} \qquad \pp$$
for any  $\gamma_1<0.35332$. 
Regarding upper bounds,  Maier \& Tenenbaum (see \citer{MT85},  \citer{MT09}) proved that, given any $\gamma_2>\log 2\approx 0.693147$, we have
$$\Delta(n)\leqslant (\log_2n)^{\gamma_2} \qquad{\rm pp.}\eqdef{majnormale}$$
We are now able to improve on this result by reducing the exponent further.
\vskip-3mm
\Propt{th2}{Let $\gamma_3>(\log 2)/(\log 2+1/\log 2-1 )\approx 0.6102495$. We have
$$\Delta(n)\leqslant (\log_2n)^{\gamma_3}\qquad {\rm pp.}$$
}
\medskip\medskip
\paraun{Average order: proof of \ref{th1}}
\paradeuxb{Reductions}
Let us start by a technical reduction similar to that of the proof of \citeplus{Te85}{(25)} and enabling to substitue the evaluation of a logarithmic mean to that of a Cesˆro mean. Considering the inequality (see \citeplus{HT88}{lemma 61.1})
$$\Delta(mn)\leqslant
\tau(m)\Delta(n)\qquad
(m\geqslant 1,\,n\geqslant 1),\eqdef{sousmultdc}$$
we may write, for any function $g$ in $\M_A(y,c,\eta)$, 
$$\sum_{n\leqslant x}g(n)\Delta(n)\log n\leqslant \sum_{m\leqslant x}g(m)\Delta(m)\sum_{\pnu\leqslant x/m}(\nu+1)g(\pnu)\log \pnu\ll x\sum_{m\leqslant x}{g(m)\Delta(m)\over m},$$
where the second bound is obtained by invoking \eqref{majgn} in the form $g(\pnu)\ll_A(3/2)^\nu$ for $p\leqslant 2A$. Since we trivially have 
$$\sum_{n\leqslant x}g(n)\Delta(n)\log \Big({x\over n}\Big)\ll x \sum_{n\leqslant x}{g(n)\Delta(n)\over n},$$
it follows that
$$\sum_{n\leqslant x}g(n) \Delta(n)
\ll{x\over \log x}\sum_{n\leqslant x}{g(n)\Delta(n)\over n}\cdot\eqdef{redlog}$$
Moreover, the canonical representation $n=md$ where $m$ is squarefree and $d$ is squarefull implies
$$\sum_{n\leqslant x}{g(n)\Delta(n)\over n}\ll \sum_{n\leqslant x}{\mu(n)^2g(n)\Delta(n)\over n}\cdot\eqdef{redsfc}$$
\par 
Next we observe that proving \eqref{majDg} for $y=1$ implies the required bound for $y\geqslant 1$.
Indeed, any $g$ in $\M_A(y,c,\eta)$ is representable as $g(n)=y^{\omega(n)}h(n)$ with
 $h\in \M_A(1,c,\eta)$.
The identity
 $$y^{\omega(n)}=\sum_{d\mid n}\mu(d)^2(y-1)^{\omega(d)},$$ already used in \citer{HT86}, and the inequality \eqref{sousmultdc}
hence imply
 $$\sum_{n\leqslant x}{g(n) \Delta(n)\over n}
 \leqslant \sum_{d\leqslant x}Ê{\mu(d)^2h(d)(2y-2)^{\omega(d)}\over d}\sum_{m\leqslant x/d} {\Delta(m)Ê h(m)\over m}\ll \e^{a\sqrt{\log_2x}}(\log x)^{2y-1},$$
by applying \eqref{majDg} to $h$. The required estimate follows by \eqref{redlog}.

In the sequel, we hence consider a  function $g\in \M_A(1,c,\eta)$ and aim at estimating the right-hand side of \eqref{redsfc}. \par \medskip
Let $\{p_j(n):1\leqslant j\leqslant \omega(n)\}$ denote the increasing sequence of distinct prime factors of a generic integer $n$ and define $n_k:=p_1(n)\cdots p_k(n)$ if $\omega(n)\geqslant k$, $n_k=n$ otherwise. Our final estimate will be obtained from a bound for 
$$D_k(x;g):=\sum_{n\leqslant x}{\mu(n)^2g(n)\Delta(n_k)\over n}$$ obtained by induction on $k$. To determine a suitable size for $k$, we appeal to a straightforward variant of \citeplus{HT88}{th. 72} providing
$$\sum_{n\leqslant x}\Delta(n)g(n)y^{\omega(n)}=x(\log x)^{2y-2+o(1)}\quad(y\geqslant 1,\,x\to\infty),$$ 
and hence, for any fixed $y>1$,
$$\sum_{\di{n\leqslant x}{\omega(n)>2y\log_2x}}g(n)\Delta(n)\ll \sum_{n\leqslant x}g(n)\Delta(n)y^{\omega(n)-2y\log_2x}\ll x(\log x)^{-2(y\log y-y+1)+o(1)}=o(x).$$
Therefore, we see that it will be sufficient to bound $D_k(x;g)$ for $k\leqslant K_x:=(2+\varepsilon)\log_2x$.
\par \medskip
A last reduction is described as follows. 
Given $\xi(x)$ tending to infinity arbitrarily slowly, let $\A_x$ denote the set of those integers $n\geqslant 1$ that are squarefree and satisfy $$\omega(n)\geqslant k\Rightarrow\log_2p_k(n)>k/5\qquad (\xi(x)\leqslant k\leqslant K_x).\eqdef{minpk}$$
Let 
$$D(x;g):=\sum_{n\leqslant x}{\mu(n)^2g(n)\Delta(n)\over n},\qquad 
D^-(x;g):=\sum_{n\in[1,x]\sset\A_x}{\mu(n)^2g(n)\Delta(n)\over n}\cdot
$$ 
\par 
Put $\omega(n,t):=\sum_{p|n,\,p\leqslant t}1$, $r_k:=\exp\exp(k/5)$.  Letting $h_{v,k}$ denote the multiplicative function supported on squarefree integers 
and defined by $h_{v,k}(p):=(v-1)\1_{[2,r_k]}(p)$, we have, for any $v\geqslant 1$
$$\eqalign{D^-(x;g) &\leqslant \sum_{\xi(x)\leqslant k\leqslant K_x}v^{-k}\sum_{n\leqslant x}{\mu(n)^2g(n)v^{\omega(n,r_k) }\Delta(n)\over n}\cr
&\leqslant  \sum_{\xi(x)\leqslant k\leqslant K_x}v^{-k}
\sum_{\di{d\leqslant x}{P^+(d)\leqslant r_k}}{\mu(d)^2g(d)h_{v,k}(d)2^{\omega(d) } \over d}
\sum_{n\leqslant x}{\mu(n)^2g(n) \Delta(n)\over n}\cr
&\ll D (x;g) \sum_{\xi(x)\leqslant k\leqslant K_x} \e^{2(v-1)k/5-k\log v}=o(D (x;g)),\cr}$$
 by selecting $v=\ft 32$ since $1/5-\log (3/2)<0$.
  \par  
 Thus we obtain that suitable averages over $\A_x$ will imply \eqref{majDg} as stated.
\medskip
\paradeuxb{A lemma}
Write
$$M_q(n):=\int_\r\Delta(n,u)^q\d u\quad(n\geqslant 1,\,q\geqslant 1).\eqdef{defMq}$$
The following estimate will play a crucial role in the proof.
\Propl{M2}{Let $A>0$, $c>0$, $\eta>0$, and $g\in\M_A(1,c,\eta)$. We have
$$\sum_{\omega(n)\geqslant k}{\mu(n)^2g(n)M_2(n_k)\over n^\sigma}\ll {  k2^k\over  \sigma-1   }\qquad (k\geqslant 1,\,1<\sigma\leqslant 2).\eqdef{majM2}$$}
\nid Let $S_2(\sigma)$ denote the left-hand side of \eqref{majM2}. Put $\tau(n,\vartheta):=\sum_{d|n}d^{i\vartheta}$ $(n\geqslant 1,\,\vartheta\in\r)$. Plainly,
$$\eqalign{S_2(\sigma)&
\ll {1\over \sigma-1}\sum_{\di{m\geqslant 1}{\omega(m)=k}}{\mu(m)^2g(m)M_2(m)\over m \log P^+(m) }
\cr&\ll {1\over \sigma-1 }\sum_{p}{g(p)\over  {p} \log p}\sum_{\di{P^+(m)<p} {\omega(m)=k-1}}{\mu(m)^2g(m) \over m } \int_\r {|\tau(m,\theta)|^2\over 1+\theta^2}\d \theta} $$
by Parseval's formula in view of \citeplus{HT88}{(3.2)}.
The last integral is classically dominated by the contribution of the interval $[-1,1]$---see the Montgomery-Wirsing lemma as stated \eg\ in \citeplus{Te15}{lemma III.4.10}.
For $|\theta|\leqslant 1$, $t\geqslant 2$, we have
$$\sum_{\di{r\leqslant t}{r\in \P}}{g(r)|\tau(r,\theta)|^2\over r}=4 \log_2t-2 \log(1+|\theta|\log t)+O(1),$$
whence
$$
\sum_{\di{P^+(m)<p} {\omega(m)=k-1}}{\mu(m)^2g(m) |\tau(m,\theta)|^2\over  m } \ll   {\{4\log_2p- 2\log( 1+|\theta|\log p)+O(1)\}^{k-1}\over    (k-1)!} \cdot$$ 
We  therefore get
$$\sum_{p}{g(p)\over  {p} \log p}\sum_{\di{P^+(m)<p} {\omega(m)=k-1}}{\mu(m)^2 g(m)\over m^\sigma } \int_\r {|\tau(m,\theta)|^2\over 1+\theta^2}\d \theta\ll{2^k\over  (k-1)!}\int_{0}^1\{T_1(\vartheta)+T_2(\vartheta)\}\d\vartheta,$$
with, for a suitable absolute constant $c_0$,
$$\eqalign{T_1(\vartheta)&:=\sum_{p\leqslant \exp(1/\vartheta)}{g(p)(2 \log_2p+c_0)^{k-1}\over p\log p},\cr
T_2(\vartheta)&:=\sum_{p>\exp(1/\vartheta)}{g(p)\{  \log_2p+ \log (1/\vartheta)+c_0\}^{k-1}\over p\log p}\cdot\cr}$$
It follows that
$$\eqalign{\int_0^1T_1(\vartheta)\d\vartheta&\ll\sum_{p}{g(p)(2 \log_2p+c_0)^{k-1}\over p(\log p)^2}\ll  (k-1)!,\cr
\int_0^1T_2(\vartheta)\d\vartheta&\ll\int_0^1
\sum_{p>\exp(1/\vartheta)}{g(p)\{ \log_2p+\log (1/\vartheta)+c_0\}^{k-1}\over p\log p}\d\vartheta\cr
&\ll  \int_0^1{1\over \vartheta}\int_{2\log (1/\vartheta)}^\infty(v+c_0)^{k-1}\e^{-v}\d v\d\vartheta\ll k!.\cr}$$
This yields \eqref{majM2} as required.\qed
\bigskip
\paradeuxb{Completion of the proof}
With  notation \eqref{defMq} and
$$ L(n):=\meas\{u\in\r:\Delta(n,u)>0\},$$ we introduce  the series 
$$
F_{k,q}(\sigma):= 
\sumast_{\omega(n)\geqslant k}{g(n)M_q(n_k)L(n_k)^{(q-1)/2}\over 2^{k } M_2(n_k)^{(q-1)/2}n^\sigma},\quad G_{k,q}(\sigma):=\sumast_{\omega(n)\geqslant k}{g(n)M_q(n_k)^{1/q}\over n^\sigma}$$
for $\sigma>1,\,k\geqslant 1$.
Here and throughout the asterisk indicates that the summation domain is restricted to $\A_x$.\par 
Since by \citeplus{HT88}{th. 72} we have $$\Delta(n)\leqslant 2M_q(n)^{1/q}\qquad (n\geqslant 1,\,q\geqslant 1),\eqdef{majDM}$$ the validity of the bound 
$$ \sum_{k\leqslant K_x} G_{k,q(k)}(1+1/\log x)\ll_a\e^{a\sqrt{\log_2x}}{\log x} \eqdef{majGkq-finale}$$ for any $a>\sqrt{2}\log 2$ and suitable $q(k)$ implies the same estimate for the right-hand side of~\eqref{redlog}.  An appropriate choice $q(k)$ will be given later.
 
Put
$$N_{j,q}(n,p):=\int_\r\Delta(n,u)^j\Delta(n,u-\log p)^{q-j}\d u,\quad W_q(n,p):=\sum_{1\leqslant j\leqslant q-1}{q\choose j}N_{j,q}(n,p).$$
The inequality
$$M_q(n_{k+1})\leqslant 2M_q(n_k)+W_q(n_k,p_{k+1})\1_{\{\omega(n)\geqslant  k+1\}},$$
is an equality if $\omega(n)\geqslant k+1$ and  holds trivially otherwise. Since $ M_2(n_{k+1})\geqslant 2M_2(n_k)$ and 
$L(n_{k+1})\leqslant 2L(n_{k})$, it follows
 that
  $$\eqalign{F_{k+1,q}(\sigma)&
\leqslant F_{k,q}(\sigma)+  \sum_{\di{m\in\M_k}{\omega(m)=k}}\sum_{\di{p>P^+(m)}{\log_2p\geqslant k/5}}{W_q(m,p)L(m)^{(q-1)/2} \over 2^kM_2(m)^{(q-1)/2}} \sumast_{n_{k+1}=mp}{g(n)\over n^\sigma},\cr}\eqdef{majdiffFk1}$$
where $\M_k:=\{m\geqslant 1:\exists n\in\A_x:n_k=m\}$.\par 
\smallskip\goodbreak
Let $\HH_k:=\{h\geqslant 1:\mu(h)^2=1,\omega(h)\geqslant j\Rightarrow \log_2p_j(h)\geqslant (j+k+1)/5\}.$
The inner sum in~\eqref{majdiffFk1} does not exceed
$${\mu(mp)^2g(mp)\over (mp)^\sigma}\sum_{\di{P^-(h)>p}{h\in\HH_k}}{g(h)\over h^\sigma},$$
hence
$$
\eqalign{&F_{k+1,q}(\sigma)- F_{k,q}(\sigma)\cr&\hskip3mm \leqslant \sum_{\di{m\in\M_k}{\omega(m)= k}}{\mu(m)^2g(m) L(m)^{(q-1)/2}\over 2^kM_2(m)^{(q-1)/2}m^\sigma} \sum_{\di{P^-(h)>P^+(m)}{h\in\HH_k}}{g(h)\over h^\sigma}
\sum_{\di{p>P^+(m)}{\log_2p\geqslant k/5}}{g(p)W_q(m,p) \over p}\cdot\cr}\eqdef{rec1}
$$
\par
Now, observe that, for all $z>1$ and $1\leqslant j<q$,
$$\eqalign{(\log z) \sum_{p>z}{g(p)N_{j,q}(m,p)\over p}&\leqslant \sum_{p>z}{g(p)N_{j,q}(m,p)\log p\over p}
\cr&=\int_\r\Delta(m,u)^j \sum_{d_1,\ldots, d_{q-j}\mid m}
\sum_{\di{p>z}{\max_h d_h<p\leqslant \e \min_h d_h}}{g(p) \log p\over p} \d u.}$$
The inner $p$-sum does not exceed $$ \Big\{1-\log \Big({\max_h d_h\over \min_h d_h }\Big) +O\Big(\e^{-c(\log
z)^\eta}\Big)\Big\}\1_{\{\max_h d_h /\min_h d_h<\e\}},$$ and so
$$\eqalign{(\log z)\sum_{p>z}{g(p)N_{j,q}(m,p)\over p}&\leqslant \sum_{p>z}{g(p)N_{j,q}(m,p)\log p\over p}
\cr&\leqslant   M_j(m)M_{q-j}(m)+O\Big(\e^{-c(\log
z)^\eta}M_j(m)M_{q-j}^*(m)\Big),}$$
where $$M_\ell^*(n):=\sum_{\di{d_h|n\ (1\leqslant h\leqslant \ell)}{\max d_h\leqslant \e\min d_h }}\leqslant 2^{\ell}M_\ell(n)\qquad (\ell\geqslant 1,\,n\geqslant 1),$$
by \citeplus{MT85}{(6)}. 
\par 
At this stage we note that H\"older's inequality furnishes for $2\leq\ell\leqslant q-2,\,m\geqslant 1$
$$M_\ell(m)\leqslant M_2(m)^{(q-\ell-2)/(q-4)}M_{q-2}(m)^{(\ell-2)/(q-4)} .$$
Applying twice  for $\ell=j$ and $\ell=q-j$, we get  for any $2\leqslant j\leqslant q-2$
$$M_j(m)M_{q-j}(m)\leqslant M_2(m)M_{q-2}(m) $$
so that  for any $2\leqslant j\leqslant q-2$
$$\sum_{p>z}{g(p)N_{j,q}(m,p)\log p\over p}\leqslant  M_2(m)M_{q-2}(m)\big\{1+R_{j,q}(m,z)\big\}$$
where
$$R_{j,q}(m,z)\ll 2^{q-j}\e^{-c(\log z)^\eta}.$$
\par 
Carrying back into \eqref{rec1} and taking into account the fact that, when $m\in\M_k$, $h\in\HH_k$, $P^-(h)>P^+(m)$, $\mu(mh)^2=1$, the integer $mh$ belongs to $\A_x$, we obtain
$$\eqalign{F_{k+1,q}(\sigma)- F_{k,q}(\sigma)\ll& q\big(1+2^q\varepsilon_k\big)H_{k,q}(\sigma)+\big(2^q+3^q\varepsilon_k\big)J_{k,q}(\sigma),\cr}$$
with
$$\eqalign{\varepsilon_k&:=\e^{-c\exp(\eta k/5)},\cr
H_{k,q}(\sigma)&:=\sumast_{\omega(n)\geqslant k}{g(n_k)M_{q-1}(n_k) L(n_k)^{(q-1)/2}\over  M_2(n_k)^{(q-1)/2}n^\sigma \log p_k(n)},\cr J_{k,q}(\sigma)&:= \sumast_{\omega(n)\geqslant k}{g(n_k)M_{q-2}(n_k) L(n_k)^{(q-1)/2}\over 2^{k }M_2(n_k)^{(q-3)/2}n^\sigma \log p_k(n)}\cdot\cr}$$

From the inequalities $${ L(n_k)^{(q-1)/2}\over   \log p_k(n)}\leqslant k   L(n_k)^{(q-3)/2},\quad{1\over L(n_k)}\leqslant {M_2(n_k)\over 4^k}  \qquad(\omega(n)\geqslant k), $$
 we deduce that
$$\eqalign{H_{k,q}(\sigma)&\leqslant \sumast_{\omega(n)\geqslant k}{kg(n_k)M_{q-1}(n_k) L(n_k)^{(q-2)/2}\over 2^{k  }M_2(n_k)^{(q-2)/2}n^\sigma }=kF_{k,q-1}(\sigma),\cr
J_{k,q}(\sigma)&\leqslant  \sumast_{\omega(n)\geqslant k}{kg(n_k)M_{q-2}(n_k) L(n_k)^{(q-3)/2}\over 2^{k  }M_2(n_k)^{(q-3)/2}n^\sigma  }=kF_{k,q-2}(\sigma),\cr}$$
whence
$$F_{k+1,q}(\sigma)- F_{k,q}(\sigma)\ll kq\big(1+2^q\varepsilon_k\big) F_{k,q-1}(\sigma)+k\big(2^{q}+3^q\varepsilon_k\big) F_{k,q-2}(\sigma). $$
 \par 
Let $k_0(q):=B\log q$, where $B$ is sufficiently large to ensure  $\varepsilon_k 2^q\leqslant 1$ whenever $q\geqslant 2$, $k\geqslant k_0(q)$. We thus have, for a suitable constant $C>0$,
$$F_{k+1,q}(\sigma)- F_{k,q}(\sigma)\leqslant Ckq F_{k,q-1}(\sigma)+ Ck2^q F_{k,q-2}(\sigma)
\qquad (k\geqslant k_0(q),\,\sigma>1).\eqdef{HR}$$ 
We now show by induction on $k$ that this implies, for a suitable constant~$D$,
$$F_{k,q}(\sigma)\leqslant {Dq^{Bq}k^{3q}2^{ q^2/4}\over \sigma-1}\prod_{1\leqslant j\leqslant k}\bigg(1 
+{4C\over j^2}\bigg) \qquad (q\geqslant 2,\,k\geqslant 1,\,1<\sigma\leqslant 2) .\eqdef{rec}$$
For $k\leqslant k_0(q)$, this follows from trivial bound
$(\sigma-1) F_{k,q}(\sigma)\ll 2^{k(q-1)}$, in view of the inequalities $1/M_2(n)\leqslant L(n)/\tau(n)^2$,  $L(n)\leqslant \tau(n)$, $ M_q(n)\leqslant  \tau(n)^q $. 
 Assuming that \eqref{rec} holds for $k\geqslant k_0(q)$, we deduce from \eqref{HR} that
  $$\eqalign{F_{k+1,q}(\sigma)
&\leqslant 
  {Dq^{Bq}(k+1)^{3q}2^{  q^2/4}\over \sigma-1} \prod_{1\leqslant j\leqslant k+1}\bigg(1 
+{4C\over j^2}\bigg), \cr}
$$ 
since $q\leqslant 2^{(q+1)/2}$ and $$k^{3q}2^{q^2/4}+Ck^{3q-2}2^{(q+1)/2+(q-1)^2/4}+ Ck^{3q-2}2^{q+(q-2)^2/4}\leqslant (k+1)^{3q}2^{q^2/4}\Big\{1+{4C\over (k+1)^2}\Big\} .$$
Therefore, we may state that
$$F_{k,q}(\sigma)\ll {q^{Bq}k^{3q}2^{ q^2/4}\over \sigma-1}  \qquad (q\geqslant 2,\,k\geqslant 1,\,1<\sigma\leqslant 2).\eqdef{majFkq}$$

 \bigskip  \par 
Now invoking once more the inequality $4^k=\tau(n_k)^2\leqslant L(n_k)M_2(n_k) $ we get via Hšlder's inequality that
$$\eqalign{G_{k,q}(\sigma)&\leqslant  2^{k/q} F_{k,q}(\sigma)^{1/q}\Bigg\{\sum_{\omega(n)\geqslant k}{\mu(n)^2g(n)\sqrt{M_2(n_k)}\over n^\sigma \sqrt{L(n_k)}}\Bigg\}^{1-1/q}\cr&\leqslant  { q^Bk^32^{k/q+q/4}\over (\sigma-1)^{1/q}}\Bigg\{\sum_{\omega(n)\geqslant k}{\mu(n)^2g(n)M_2(n_k)\over 2^{k}n^\sigma}\Bigg\}^{1-1/q}\ll{q^Bk^42^{k/q+q/4}\over \sigma-1},\cr}\eqdef{majG}$$
by \eqref{majM2}.
\par

Applying this for all $k\leqslant K_x$, with $q=q(k)=\lceil 2\sqrt{ k}\,\rceil$ and $\sigma=1+1/\log x$, yields \eqref{majGkq-finale} and thus finishes the proof.
\goodbreak
\paraun{Normal order: proof of \ref{th2}}

This is a reappraisal of \citeplus{MT09}{th. 1.3}. By \eqref{majDM}, it is sufficient to bound  $M_q(n)$.
\par 
Let $\lambda\in ]1,2[$ and let $\gamma$, $\delta$ be real numbers satisfying
$$ \delta(\log 2)/\lambda<\gamma<1 ,\qquad 1 <\delta < \lambda(\gamma-1) +1/\log 2.\eqdef{encdelta}$$
We shall show by induction on $k$ that 
$$M_q(n_k)\leqslant 2^{\delta k} (q!)^\gamma\qquad (1\leqslant q\leqslant \lambda k) \quad \pp x.\eqdef{hr}$$
Here, as in \citer{HT88}, the mention $\pp x$ indicates that a formula holds for all but at most $o(x)$ integers $n\leqslant x$ as $x\to\infty$.\par 
Given an integer-valued function
$\xi=\xi(x)$ tending to
infinity arbitrarily slowly, we put 
$$K=K(n,x):=\max\{k:1\leqslant  k\leqslant  \omega(n),\,\log_2p_k(n)<\log_2x-\xi(x)\} $$
and redefine
$$n_k:=\cases{\prod_{\xi<j\leqslant  k}p_j(n)& if $k\leqslant  K,$\cr n_K & if $k>K$.\cr}$$
\par 
Assuming \eqref{hr} holds for  $k$ with  $\xi<k\leqslant K$ we aim at showing that this bound persists at rank $k+1$.
\par 
Let $\e_1<\e$. By \citeplus{MT09}{(3.2)}, for $\xi<k\leqslant K,\, q\geqslant 1$, we have
$$M_q(n_{k+1})  \leqslant 2M_q(n_{k})+\e_1^{-k}
\sum_{1\leqslant j\leqslant q-1} {q\choose j}M_j(n_{k})M_{q-j}(n_{k})\quad {\rm pp}x.\eqdef{3.2}$$
Note that H\"older's inequality implies 
$$\sum_{1\leqslant j\leqslant q-1} {q\choose j}^{1-\gamma}\leqslant 2^{(1-\gamma)q}(q-1)^\gamma.$$
\par 
When $q\leqslant q_k:=\fl{\lambda k}$,  we appeal to the induction bound \eqref{hr} to majorize the right-hand side of \eqref{3.2}. This yields
$$\eqalign{M_q(n_{k+1}) 
&\leqslant 2^{\delta (k+1)} (q!)^\gamma\bigg\{2^{1-\delta}+(2^\delta/\e_1)^{ k}
\sum_{1\leqslant j\leqslant q-1} {q\choose j}^{1-\gamma}\bigg\}
\cr&\leqslant 2^{\delta (k+1)} (q!)^\gamma\Big\{2^{1-\delta}+(2^\delta/\e_1)^{ k}
2^{(1-\gamma)q}q^\gamma\Big\}
\cr&\leqslant 2^{\delta (k+1)} (q!)^\gamma\Big\{2^{1-\delta}+(2^{\delta+(1-\gamma)\lambda}/\e_1)^{ k}
q^\gamma\Big\}\leqslant 2^{\delta (k+1)} (q!)^\gamma,
\cr}\eqdef{qqk}$$
for sufficiently large $\xi$, by the second condition \eqref{encdelta}.  
\goodbreak
When $q_k <q\leq\fl{\lambda (k +1)}$,  we apply  \citeplus{MT09}{(3.3)}: given $\alpha>0$ and $r>1/\alpha$ we have
$$\Delta(n_k)\leqslant r+\e^{\alpha k/q} M_q(n_k)^{1/q} \qquad (\xi<k\leqslant K, q\geqslant 1)\quad\pp x.$$

Since $\lambda<2$, we have $q=q_k+1$ or $q=q_k+2$. If $r$ is sufficiently large and $\alpha >1/r$, our induction hypothesis \eqref{hr} furnishes
$$\eqalign{M_{q_k+1}(n_{k}) &\leqslant  \Delta(n_k)M_{q_k}(n_{k})
\leqslant 
\big\{ \e^{\alpha/\lambda}M_{q_k}(n_{k})^{1/q_k}+r\big\} M_{q_k}(n_{k})
\cr&
\leqslant 2^{\delta k} \{(q_k+1)!\}^\gamma\bigg\{{ 2^{ \delta/\lambda }\e^{\alpha/\lambda}(q_k!)^{\gamma/q_k }\over (q_k+1)^\gamma } +{r\over (q_k+1)^\gamma } \bigg\}
\cr&
\leqslant 2^{\delta k} \{(q_k+1)!\}^\gamma\big\{2^{ \delta/\lambda }\e^{\alpha/\lambda-\gamma}+o(1)\big\} .\cr}
\eqdef{q=2k+1,k}$$
 
Carrying back into  \eqref{3.2} and writing $$b:=2^{1-\delta+\delta/\lambda}\e^{\alpha/\lambda-c}<1,$$ 
we get, taking \eqref{qqk} into account,
$$\eqalign{M_{q_k+1}(n_{k+1}) & \leqslant 2M_{q_k+1}(n_{k})+\e_1^{-k}
\sum_{1\leqslant  j\leqslant q_k}{q_k+1 \choose j}M_j(n_{k})M_{{q_k+1}-j}(n_{k})
\cr&\leqslant  2^{\delta (k+1 )} \{(q_k+1)!\}^\gamma\bigg\{ b+o(1)+(2^\delta/\e_1)^{ k}
\sum_{1\leqslant j\leqslant q_k} {q_k+1\choose j}^{1-\gamma}\bigg\}
\cr&\leqslant 
2^{\delta (k+1 )} \{(q_k+1)!\}^\gamma\bigg\{b+o(1)+(2^\delta/\e_1)^{ k}
2^{(1-\gamma)( \lambda k  +1)}(q_k+1)^\gamma\bigg\}
\cr&\leqslant   
2^{\delta (k+1)} \{(q_k+1)!\}^\gamma\big\{b+o(1)\big\}
\leqslant 2^{\delta (k+1)} \{(q_k+1)!\}^\gamma,\cr}\eqdef{qqk+1}$$
by \eqref{encdelta}.
\par 
If $q_k+2\leqslant  \lambda (k+1)$, we also need to bound $M_{q_k+2} (n_{k+1})$. Put $g:=2^{\delta/\lambda}\e^{\alpha/\lambda-c}<1.$ By \eqref{hr}, \citeplus{MT09}{(3.3)} and \eqref{qqk+1}, we have, for large $\xi$,
$$\eqalign{
 M_{q_k+2} (n_{k+1}) &\leqslant \Delta (n_{k+1})M_{q_k+1}(n_{k+1})\cr&\leqslant 
\big\{r+ \e^{\alpha/\lambda}M_{q_k+1}(n_{k+1})^{1/(q_k+1)}\big\} M_{q_k+1}(n_{k+1})
\cr&
\leqslant 2^{\delta (k+1)} \{(q_k+2)!\}^\gamma{b+o(1)\over (q_k+2)^\gamma}\bigg\{r+\e^{\alpha/\lambda}2^{\delta/\lambda+\delta/\lambda k}\{(q_k+1)!\}^{\gamma/(q_k+1) } \bigg\}
\cr&
\leqslant 2^{\delta (k+1)} \{(q_k+2)!\}^\gamma\{b+o(1)\}\big\{g+o(1)\big\}
\leqslant  2^{\delta (k+1)} \{(q_k+2)!\}^\gamma,\cr}
 $$
still by \eqref{encdelta}.  
\par \goodbreak
Selecting   $\alpha $ sufficiently small, we see that the induction hypothesis is still valid at rank $k+1$.  
We may take $\delta $ arbitrarily close to  $1$, and so $\gamma$ arbitrarily close to $\gamma_3$. This yields 
$$M_{q_K}(n_K)\leqslant 2^{\delta K/q_K} (q_K!)^{\gamma/q_K}\ll K^{\gamma}.$$
Since we have classically $K(n,x)\sim \log_2x$ $\pp x$ as $x\to\infty$,
we may conclude as in \citer{MT09} by invoking the bound
$$\Delta(n)\leqslant  \Delta(n_K)2^{\Omega(n/n_K)} \ll \Delta(n_K)4^{\xi}\qquad \pp x.$$
 \bigskip

 \centerline{\twelvebf References}
\medskip
 {\eightpoint\leftskip9mm\rightskip5mm
 \bibtem{Er73} P. Erd\H os, Problem 218, Solution by the proposer, {\it Can. Math. Bull. \bf 16} (1973),  621-622.\par 
 \bibtem{EN75} P. Erd\H os \& J.-L. Nicolas, RŽpartition des nombres superabondants, {\it Bull. Soc. math. France \bf 103} (1975), 65-90.\par 
 \bibtem{FGK19} K. Ford,  B. Green, D. Koukoulopoulos, Equal sums in random sets and the concentration of divisors, preprint, 2019.
 \bibtem{HT86} R.R. Hall \& G. Tenenbaum, The average orders of Hooley's $\Delta_r$-functions, II,
 {\it Compositio Math. \bf 60} (1986), 163--186.\par 
 \bibtem{HT88} R.R. Hall \& G. Tenenbaum, {\it Divisors}, Cambridge tracts in
mathematics, no 90, Cambridge University Press (1988).\par 
\bibtem{Ho79} C. Hooley, On a new technique and its applications to the theory of numbers, {\it Proc.
London Math. Soc.} (3) {\bf 38} (1979), 115--151.
\par
\bibtem{MT84}  H. Maier \& G. Tenenbaum,  On the set of divisors of an integer, {\it Invent. Math.} {\bf 76}
(1984), 121--128. 
 \bibtem{MT85}  H. Maier \& G. Tenenbaum, On the normal concentration of divisors, {\it J. London Math. Soc.} (2) {\bf 31} (1985), 393--400. \par 
  \bibtem{MT09} H. Maier and G. Tenenbaum, On the normal concentration of divisors, 2, {\it Math. Proc. Camb. Phil. Soc. \bf147} \numero3 (2009), 593-614.\par
 \bibtem{Ro11}   O. Robert, 
Sur le nombre des entiers reprŽsentables comme somme de trois puissances, {\it Acta Arith. } {\bf 149} (2011), \numero1, 1--21. 
 \bibtem{Sh80} P. Shiu, A Brun--Titchmarsh theorem for multiplicative functions, {\it J. reine
angew. Math. \bf313} (1980), 161--170.
\bibtem{Te85} G. Tenenbaum,
Sur la concentration moyenne des diviseurs, {\it Comment. Math. Helvetici} {\bf 60} (1985), 411-428. 
\bibtem{Te86} G. Tenenbaum,
Fonctions $\Delta$ de Hooley et applications, SŽminaire de thŽorie des nombres, Paris 1984-85, Prog. Math. 63 (1986), 225-239.
 \par
\bibtem{Te15} G. Tenenbaum, {\it Introduction to analytic and probabilistic number theory}, 3rd ed., Graduate Studies in Mathematics 163, Amer. Math. Soc. (2015).\par}
{\leftskip2mm\rightskip-2cm\sevenrm
\gutter=4cm \doublecolumns
 \obeylines \baselineskip=7pt
RŽgis de la Bretche
UniversitŽ Paris CitŽ, Sorbonne UniversitŽ, CNRS 
Institut de Math. de Jussieu-Paris Rive Gauche
 F-75013 Paris  
France
\smallskip
{\seventt regis.de-la-breteche@imj-prg.fr}
\phantom a\goodbreak
 GŽrald Tenenbaum
Institut \'Elie Cartan
Universit\'e de Lorraine
 BP 70239\par
54506 Vand\oe uvre-ls-Nancy Cedex
 France
\smallskip
{\seventt gerald.tenenbaum@univ-lorraine.fr}
\singlecolumn
\par}

\bye